\documentclass[preprint]{elsarticle}

\usepackage{amsmath,amssymb}
\usepackage{amsthm}
\usepackage{latexsym}
\usepackage{graphicx}
\usepackage{color}
\usepackage{subfigure}
\usepackage{hyperref}
\usepackage{url}

\newcommand{\NN}{\mathbb{N}}

\newcommand{\RR}{\mathbb{R}}

\newcommand{\EDR}{\textup{EDR}}

\newcommand{\ECG}{\textup{ECG}}
\newcommand{\II}{\textup{II}}

\theoremstyle{definition}

\theoremstyle{remark}

\begin{document}
\begin{frontmatter}
\title{Using synchrosqueezing transform to discover breathing dynamics from ECG signals}

\author[HTWu]{Hau-tieng Wu}
\ead{hauwu@stanford.edu}

\author[YHChan]{Yi-Hsin Chan}

\author[YTLin]{Yu-Ting Lin}

\author[YHYeh]{Yun-Hsin Yeh}

\address[HTWu]{Department of Mathematics, Stanford University, Stanford, USA.}
\address[YHChan]{First Division of Cardiovascular Department, Chang Gung Memorial Hospital, Linkou, Taoyuan, Taiwan.}
\address[YTLin]{Department of Anesthesiology, Shin Kong Wu Ho-Su Memorial Hospital, Taiwan and Graduate Institute of Biomedical Electronics and Bioinformatics, National Taiwan University, Taiwan}
\address[YHYeh]{First Division of Cardiovascular Department, Chang Gung Memorial Hospital, Linkou, Taoyuan, Taiwan and College of Medicine, Chang Gung University, Taoyuan, Taiwan.}

\begin{abstract}
The acquisition of breathing dynamics without directly recording the respiratory signals is beneficial in many clinical settings. The electrocardiography (ECG)-derived respiration (EDR) algorithm enables data acquisition in this manner. However, the EDR algorithm fails in analyzing such data for patients with atrial fibrillation (AF) because of their highly irregular heart rates. To resolve these problems, we introduce a new algorithm, referred to as SSTEDR, to extract the breathing dynamics directly from the single lead ECG signal; it is based on the EDR algorithm and the time-frequency representation technique referred to as the synchrosqueezing transform. We report a preliminary result about the relationship between the anesthetic depth and breathing dynamics. To the best of our knowledge, this is the first algorithm allowing us to extract the breathing dynamics of patients with obvious AF from the single lead ECG signal.
\end{abstract}

\end{frontmatter}

\section{Introduction}
Respiratory signals can contain a wealth of information valuable for the clinics. For example, the existence of the respiratory rate variability is commonly observed in normal subjects, and its clinical application in ventilator weaning prediction has been explored extensively recently \cite{Wysocki2006}.    
We shall call the information extracted from the respiratory signal of this kind {\it breathing dynamics}.

However, obtaining the respiratory signal might be a cumbersome task,  
and may not be convenient for certain clinical purposes, such as ambulatory, monitoring of anesthesia for surgical procedures, or long-term monitoring in naturalistic settings.  
Based on the profound information contained in the respiratory signal, the development of a convenient method to obtain such signals or their dynamics, is important from a clinical perspective.

\begin{figure}[!t]
\centering
\subfigure{
\includegraphics[width= \textwidth]{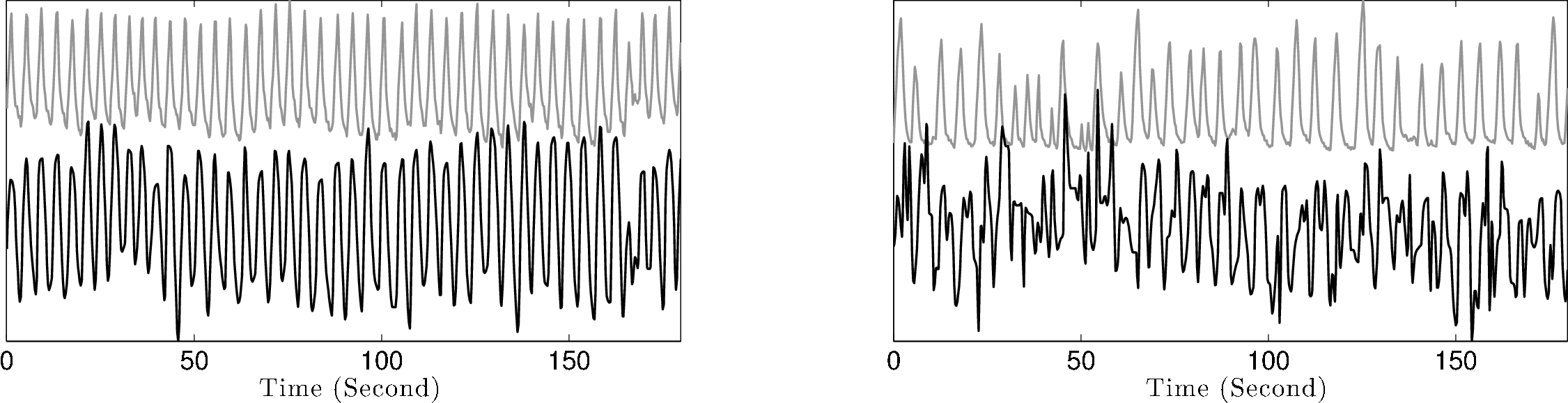}}
\caption{Left: The EDR signal, $Y_{\textup{EDR}}$, of a subject with normal sinus rhythm is illustrated by the black curve; the respiration signal measured from the chest belt, $Y_{\textup{THO}}$, is illustrated by the gray curve. Note that the $Y_{\textup{EDR}}$ correlates well with the $Y_{\textup{THO}}$. Right: The $Y_{\textup{EDR}}$ of a subject with AF is illustrated by the black curve; the $Y_{\textup{THO}}$ is illustrated by the gray curve. 
Clearly, the $Y_{\textup{EDR}}$ corresponds poorly with the $Y_{\textup{THO}}$. }
\label{fig:1}
\end{figure}

The ECG-derived respiration ($\EDR$) algorithm is one way to achieve this goal, in that it allows us to obtain a respiratory signal from the easily accessed ubiquitous ECG signals.  EDR algorithm is based on respiration induced ECG distortions, which are caused by $2$ physiological reasons. 
First, the respiration-related mechanical changes affect the thoracic electrical impedance \cite{Clifford_Azuaje_McSharry:2006}, and the cardiac axis rotation occurring during the respiratory cycle has been shown to be the largest factor contributing to the distortion of ECG signals \cite{moody_mark:1986}. Secondly, respiration affects the heart rate variability (HRV), thereby causing respiratory sinus arrhythmia \cite{Clifford_Azuaje_McSharry:2006}. 
Many $\EDR$ algorithms have been developed to obtain the respiratory signal on the basis of the 2 physiological factors mentioned above, and hence are essentially classified into $2$ major groups.  
With the acquired EDR signals, we can extract the breathing dynamics of the subject.  
However, as far as we know, the existing EDR algorithms fail in patents with atrial fibrillation (AF) due to  
their white-noise-like HRV (see example in Figure \ref{fig:1}).  
Thus, it may not be possible to obtain breathing dynamics from the EDR signals in AF patients. In light of the high prevalence of AF in clinics (it affects more than $5\%$ of individuals older than $65$ years) and the ubiquitousness of ECG, it is beneficial to develop a new algorithm that allows us to extract the breathing dynamics directly from the ECG signal, which is particularly suitable for AF patients.

We thus propose a novel algorithm for this purpose, which comprises the EDR algorithm and an adaptive time-frequency analysis technique, referred to as synchrosqueezing transform (SST). We shall refer the algorithm as SSTEDR.  
To demonstrate the applicability of SSTEDR, the algorithm is applied to 
reveal the relationship between the anesthetic depth and breathing dynamics.

\section{Mathematical Model and Algorithm}\label{section:model} 
\subsection{Mathematical Model}

We model the {\em measured respiratory signal} by considering the equation
\begin{equation}\label{eq:resp_model}
Y_{\{\}}(t)=A(t)s_{\{\}}(\phi(t))+\Phi_{\{\}}(t),
\end{equation}
where ${\{\}}$ indicates the different respiratory signal measurements, for example, $Y_{\textup{Cflow}}$ represents the airflow measuring device, $Y_{\textup{THO}}$ represents the chest wall movement measuring device, $Y_{\textup{EDR}}$ represents the traditional EDR signal, etc, and the following conditions are satisfied:
\begin{enumerate}
\item[(I)] $A\in C^1(\RR)\cap L^\infty(\RR)$, $\phi\in C^2(\RR)$, $\inf_{t\in\RR} A(t)>0$, $\inf_{t\in\RR}\phi'(t)>0$, $|A'(t)|\leq \epsilon \phi'(t)$ and $|\phi''(t)|\leq \epsilon\phi'(t)$, for all $t$, and $\epsilon$ is a small parameter;
\item[(II)] $s_{\{\}}:[0,1]\to \RR$ is $C^{1,\alpha_{\{\}}}$, where $\alpha_{\{\}}>1/2$, and $1$-periodic function with unit $L^2$ norm, $|\widehat{s_{\{\}}}(k)|\leq \delta_{\{\}} |\widehat{s_{\{\}}}(1)|$ for all $k\neq 1$, where $\delta_{\{\}}\geq 0$ is a small parameter, and $\sum_{n>D_{\{\}}}|n\widehat{s_{\{\}}}(n)|\leq \theta_{\{\}}$ for some small parameter $\theta_{\{\}}\geq 0$ and $D_{\{\}}\in \NN$. Here $\alpha_{\{\}}$, $\delta_{\{\}}$, $\theta_{\{\}}$ and $D_{\{\}}$ depend on the measurement device;
\item[(III)] $\Phi_{\{\}}(t)$ is a stationary random process or ``almost'' stationary random process \cite{Chen_Cheng_Wu:2013} independent of $A(t)s_{\{\}}(\phi(t))$, which is introduced to model the measurement error.  
\end{enumerate}
Note that the signal is assumed to have just $1$ component, unlike the case in previous studies wherein the signals were considered to include several components \cite{Wu:2013}. 
We shall call $s_{\{\}}(\cdot)$ the {\em wave shape function} \cite{Wu:2013}, $\phi$ the {\it phase function} of the signal $R(t)$, the derivative $\phi'(t)$ of the phase function the {\em instantaneous respiratory rate} (IRR), and $A(t)$ the {\it amplitude modulation} function (AM) \cite{Chen_Cheng_Wu:2013}.  

Physiologically, $\phi'(t)$, $A(t)$ and $s_{\{\}}$ quantify the breathing dynamics. Indeed, the breathing rate variability can be characterized by IRR, the breathing depth variance can be captured by AM, and the wave shape function models a respiratory cycle which is comprised of inspiration and an expiration.  

\subsection{The EDR Algorithm}\label{section:algorithm:EDR}

Denote the recorded lead II ECG signal as $\ECG_\II(t)$.  
We apply the median filter to remove the wandering baseline artifact in $\ECG_\II(t)$, which may come from patient movement, dirty lead electrodes and a variety of other things. We choose the moving window of length $120$ ms, so that it is longer than the average length of the QRS complex. We use the same notation $\ECG_\II(t)$ to denote the lead II ECG signal with the wandering baseline removed. 

To construct the traditional EDR signal $Y_{\EDR}(t)$, we determine the timing of the R peaks (or S-peaks when the cardiac axis is deviated) by the standard algorithm \cite{Clifford_Azuaje_McSharry:2006}. Suppose that there are $n$ normal R-peaks and their timings are $t_i$, where $i=1,\ldots,n$. Then, $Y_{\EDR}(t)$ is constructed by the cubic spline interpolation of the data points $\{(t_i,\,\ECG_\II(t_i))\}_{i=1}^n$ \cite{Clifford_Azuaje_McSharry:2006}. 
See Figure \ref{fig:1} for example .

\subsection{Synchrosqueezing Transform}\label{section:algorithm:SST}
Our main focus is to extract the breathing dynamics from the traditional EDR signal $Y_{\textup{EDR}}$, even when the subject has AF. 
The main mathematical tool we apply to achieve this goal is an adaptive time-frequency (TF) analysis, referred to as the Synchrosqueezing transform (SST) \cite{daubechies_maes:1996,Daubechies_Lu_Wu:2010,brevdo_fuckar_thakur_wu:2012,Chen_Cheng_Wu:2013}, which is a special reallocation technique \cite{Flandrin:99,Auger_Flandrin_Lin_McLaughlin_Meignen_Oberlin_Wu:2013}. We refer the reader to \cite{daubechies_maes:1996} for the initial proposal of SST. For the theoretical results, the reader can find the first theoretical analysis of SST when the signal is not contaminated by noise in \cite{Daubechies_Lu_Wu:2010} and when the signal is contaminated by stationary or almost stationary generalized random process and there exists trend in \cite{brevdo_fuckar_thakur_wu:2012,Chen_Cheng_Wu:2013}. The numerical implementation details of SST can be found in \cite{brevdo_fuckar_thakur_wu:2012}.

We summarize the theoretical results of SST relevant to this work \cite{Daubechies_Lu_Wu:2010,brevdo_fuckar_thakur_wu:2012,Chen_Cheng_Wu:2013} as follows:
\begin{itemize}
\item[(P1)] SST is robust to the several different kinds of noise $\Phi_{\{\}}$, which might be slightly non-stationary, and the estimation of IRR and AM is not influenced by the non-harmonic shape function $s_{\{\}}$;  
\item[(P2)] Since SST is local in nature, we are able to detect components that do not exist all the time and hence the dynamical behavior of the signal;  
\item[(P3)] The time-frequency representation determined by SST is visually informative; 
\item[(P4)] SST is ``adaptive'' to the data in the sense that the error in the estimation depends only on the first three moments of the mother wavelet instead of the profile of the mother wavelet.  
\end{itemize}

\subsection{The SSTEDR Algorithm}

\begin{figure}[!t]
\begin{center}
\subfigure{
\includegraphics[width= \textwidth]{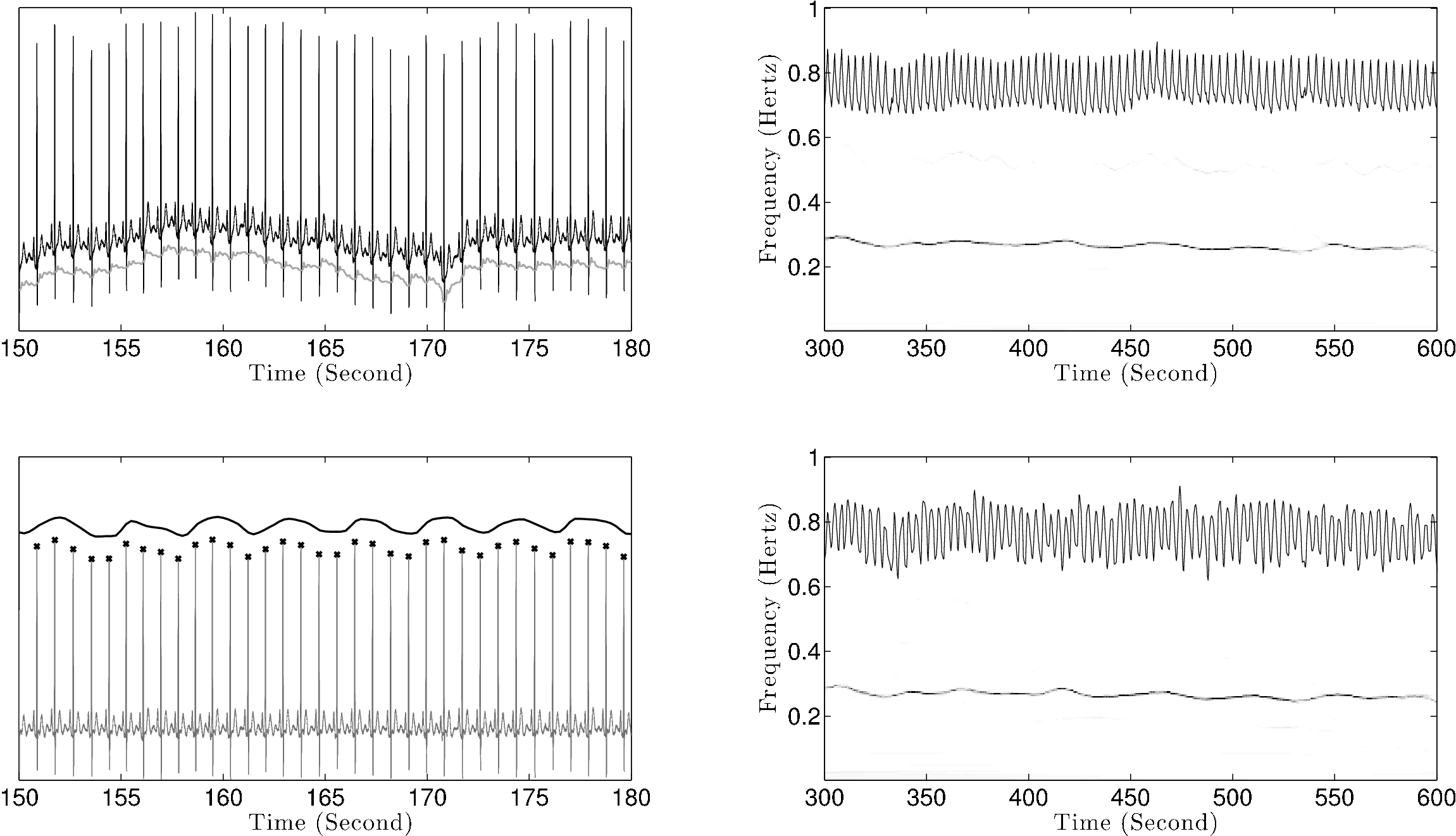}}
\end{center}
\caption{The demonstration of the SSTEDR algorithm applied to a subject with normal sinus rhythm. Upper Left: the lead II ECG signal is shown as a black curve and the wandering baseline of the lead II ECG signal determined by the median filter is plotted as a light gray curve shifted below to enhance the visualization. Lower Left: the median-filtered lead II ECG signal is plotted as a dark gray curve superimposed with the R peaks marked as black crosses, and the EDR signal, $Y_{\textup{EDR}}$, is plotted as a black curve shifted up to increase the visualization. Note that we can find the oscillatory pattern in the $Y_{\textup{EDR}}$. Right Upper: the SST of $Y_{\textup{THO}}$ with $Y_{\textup{THO}}$ superimposed as the black curve. Right Lower: the SST of $Y_{\textup{EDR}}$ with $Y_{\textup{EDR}}$ superimposed as the black curve. Notice that we can see a dominant curve inside the band between $0.2$ and $0.25$ Hertz, which is corresponding to the IRR, in the SST of $Y_{\textup{THO}}$, as well as in the SST of $Y_{\textup{EDR}}$. Furthermore, note that these two dominant curves coincide with each other.}
\label{fig:22}
\end{figure}

\begin{figure}[!t]
\begin{center}
\subfigure{
\includegraphics[width= \textwidth]{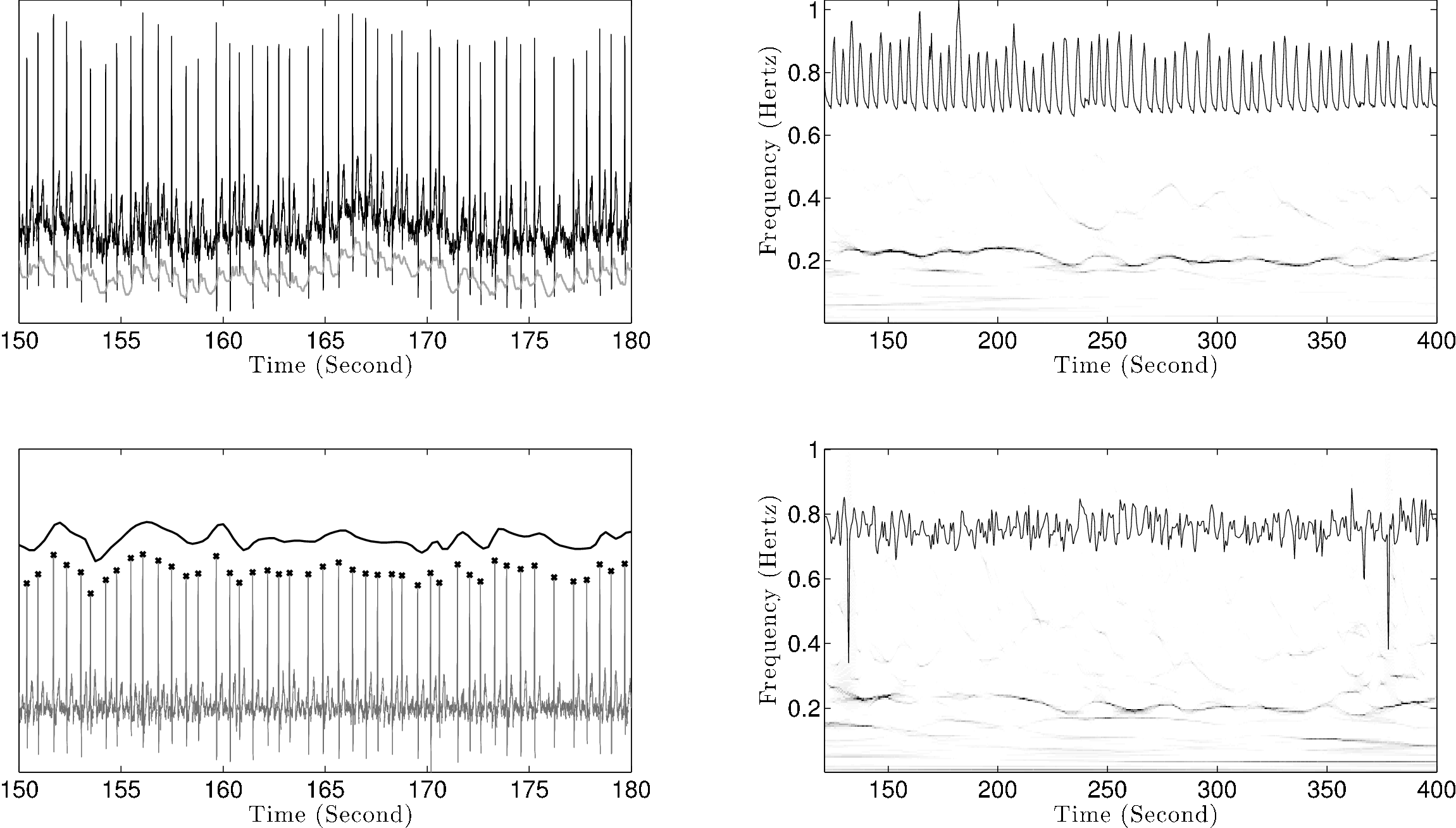}}
\end{center}
\caption{The demonstration of the SSTEDR algorithm applied to a subject with AF. Upper Left: the lead II ECG signal is shown as a black curve and the wandering baseline of the lead II ECG signal determined by the median filter is plotted as a light gray curve shifted below to enhance the visualization. Lower Left: the median-filtered lead II ECG signal is plotted as a dark gray curve superimposed with the R peaks marked as black crosses, and the EDR signal, $Y_{\textup{EDR}}$, is plotted as a black curve shifted up to increase the visualization. Note the highly irregular intervals between two consecutive R peaks. Also note that it is not easy to identify the oscillatory patterns in $Y_{\textup{EDR}}$. Right Upper: the SST of $Y_{\textup{THO}}$ with $Y_{\textup{THO}}$ superimposed as the black curve. Right Lower: the SST of $Y_{\textup{EDR}}$ with $Y_{\textup{EDR}}$ superimposed as the black curve. Notice that although $Y_{\textup{EDR}}$ and $Y_{\textup{THO}}$ have different appearance, their TF representations determined by SST contain similar information, in the sense that the dominant curves in both TF representations coincide with each other.}
\label{fig:33}
\end{figure}

The SSTEDR algorithm we propose to extract the breathing dynamics from the single lead ECG signal comprises of the following 3 steps: (1) pre-processing the lead II ECG signal to remove the wandering baseline; (2) constructing $Y_{\textup{EDR}}$ by cubic spline interpolation of the detected R peaks as is discussed in Section \ref{section:algorithm:EDR}; (3) analyze $Y_{\textup{EDR}}$ by SST. The result of SSTEDR applied to a subject without AF (resp. with AF) is illustrated in Figure \ref{fig:22} (resp. Figure \ref{fig:33}).
A Matlab code for SSTEDR can be found in \url{https://sites.google.com/site/hautiengwu/}.

\section{Anesthestic Depth Estimation}

Anesthesia is usually necessary for a patient receiving surgery. We now demonstrate a potential clinical application of SSTEDR in the anesthesia depth estimation problem. 
The dosage of anesthetic medication should be dynamically adjusted to achieve adequate level of anesthesia during surgery. In anesthesia, the respiratory pattern reflects the effects of anesthetic agents on human body. Indeed, more than seventy years ago, Guedel described that respiratory patterns appear regular in deeper levels and irregular in lighter levels of ether anesthesia \cite{Guedel:37}. Also, the respiratory pattern changes correspondingly with electroencephalography during arousal from anesthesia \cite{Bimer_Bellville:1977}. Since the clinician relies on various information to judge the administration of anesthetics, the knowledge of the breathing dynamics is appealing. 

However, while it is currently mandatory to monitor the ECG signal during anesthesia, the respiratory signal is not routinely monitored. Hence it is inconvenient to obtain the breathing dynamics during anesthesia. We therefore ask if we can apply our algorithm on the ECG signal to capture the reported breathing dynamics in \cite{Guedel:37}: the respiration is regular in deeper level of anesthesia and irregular in lighter level of anesthesia. 

The study was approved by the Institutional Review Board (Shin Kong Wu Ho-Su Memorial Hospital) and written informed consents were obtained from the subjects.  
The raw limb lead II ECG signal during anesthesia was recorded at the sampling rate $1000$ Hz and with $12$ bits resolution for off-line analysis (My ECG E3-80, MSI, Inc., Taiwan).    
The inhaled and end-tidal concentrations of the anaesthetic gas detected by the gas analyser on a Datex-Ohmeda S/5 anaesthesia machine (GE Health Care, Helsinki, Finland) were simultaneously recorded. The airflow signal was continuously measured by the anesthesia machine too. Data from the anesthesia machine were recorded at a rate of $25$Hz on the laptop via a dedicated software (Datex-Ohmeda S/5 collect. Ver 4, GE Health Care, Helsinki, Finland).

\begin{figure}[!t]
\begin{center}
\subfigure{
\includegraphics[width= \textwidth]{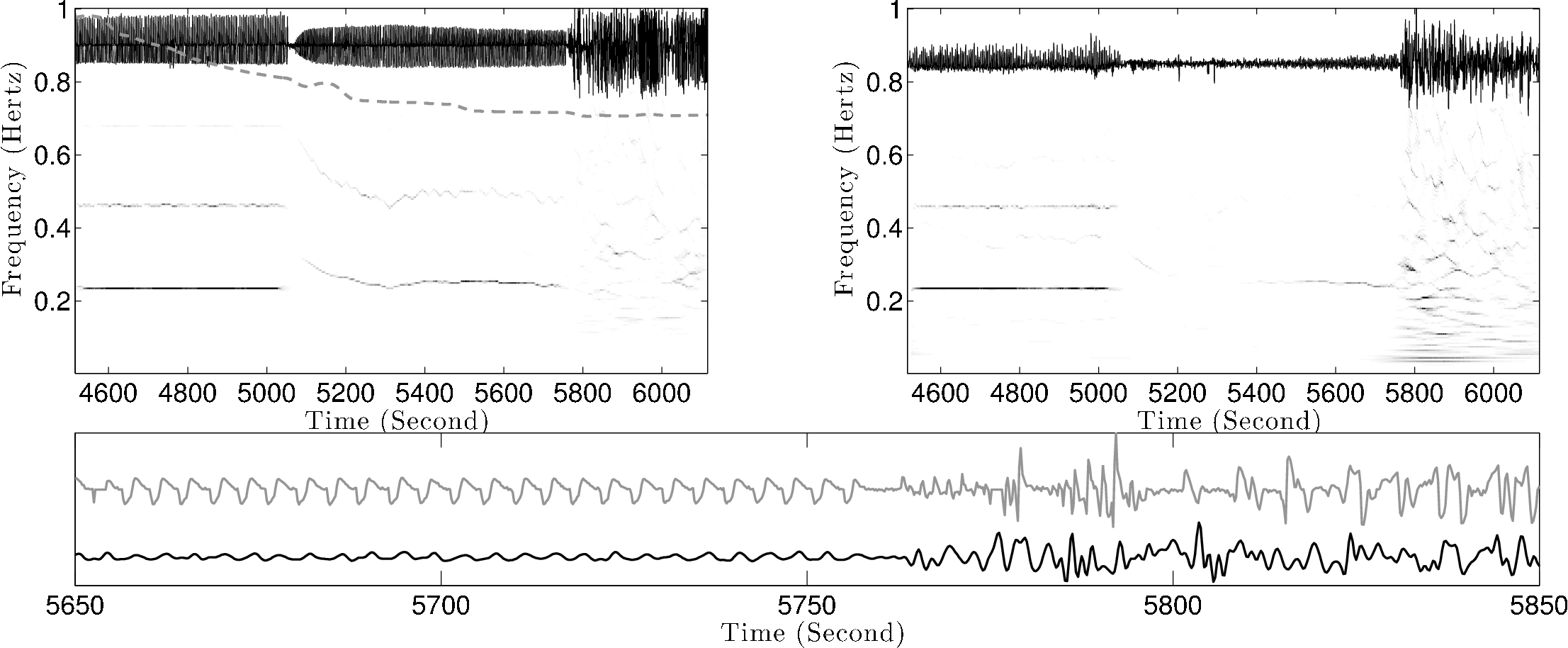}} 
\end{center}
\caption{Upper Left: the SST of the recorded respiratory signal, $Y_{\textup{Cflow}}$, from a subject with normal sinus rhythm with $Y_{\textup{Cflow}}$ superimposed as the black curve and the end-tidal concentrations of anaesthetic gas superimposed as the gray dashed curve. The main surgical procedure was finished at the 4500-th second and the skin suture was started. In the mean time, the anesthesiologist started to reduce the anesthetic medication concentration, as can be seen in the red curve.   
The mechanical ventilator was removed at the 5050-th second and the subject started to breath spontaneously. At the 5781-th second, the first reaction (movement of his head and neck) occurred.  Notice the dominant curve inside the band around $0.25$ Hertz before the $5781$-th second is corresponding to the IRR. The signal after the first reaction, however, does not satisfy the model (\ref{eq:resp_model}), so we can not see any dominant curve and the TF representation of SST is blurred. Upper Right: the SST of $Y_{\textup{EDR}}$, superimposed with $Y_{\textup{EDR}}$, as the black curve. Note that the pattern transition can be seen both in the SST of $Y_{\textup{EDR}}$, as well as the waveform of $Y_{\textup{EDR}}$. Also note that the dominant curve existing from the 5100-th second to the 5300-th second in the SST of $Y_{\textup{Cflow}}$ is not that clear in the SST of $Y_{\textup{EDR}}$. It is because the amplitude of $Y_{\textup{EDR}}$, is relatively small during this period compared with other periods.
Lower: the zoomed in signal around the first reaction (at the 5781-th second) of a subject with normal sinus rhythm. The $Y_{\textup{Cflow}}$ is shown in the gray curve and the $Y_{\textup{EDR}}$ is shown in the black curve. Note that we can easily see visually the transition from ``regular oscillation'' to ''irregular oscillation'' in both signals. The regular oscillation is modeled by (\ref{eq:resp_model}).}
\label{fig:10}
\end{figure}

We show the analysis result of an anesthetized patient with normal sinus rhythm under anesthesia in Figure \ref{fig:10}.
The SST of the recorded respiratory flow $Y_{\textup{Cflow}}$ is shown in the upper left subfigure of Figure \ref{fig:10} with $Y_{\textup{Cflow}}$ superimposed as the black curve and the end-tidal concentrations of the anaesthetic gas superimposed as the red curve.   
As the concentration of anesthetic gas decreased, gradually the subject awakened from anesthesia and the motor movement was recovered. 
The statement that the respiration is regular in the deeper anesthesia level, that is, before the first reaction, and irregular in the lighter anesthesia level, that is, after the first reaction, can be visually seen in the SST of $Y_{\textup{Cflow}}$.    
Indeed, we can see a dominant curve in the TF representation before the first reaction, while the TF representation becomes blurred afterward. Note that this transition from the regular pattern to non-regular pattern can be seen clearly in the recorded respiratory signal, too, as is shown in the lower subfigure of Figure \ref{fig:10}. On the other hand, we can also visually observe this pattern transition in $Y_{\textup{EDR}}$ shown in the lower subfigure of Figure \ref{fig:10} as well as its SST shown in the upper right subfigure of Figure \ref{fig:10}.

\begin{figure}[!t]
\begin{center}
\subfigure{
\includegraphics[width= \textwidth]{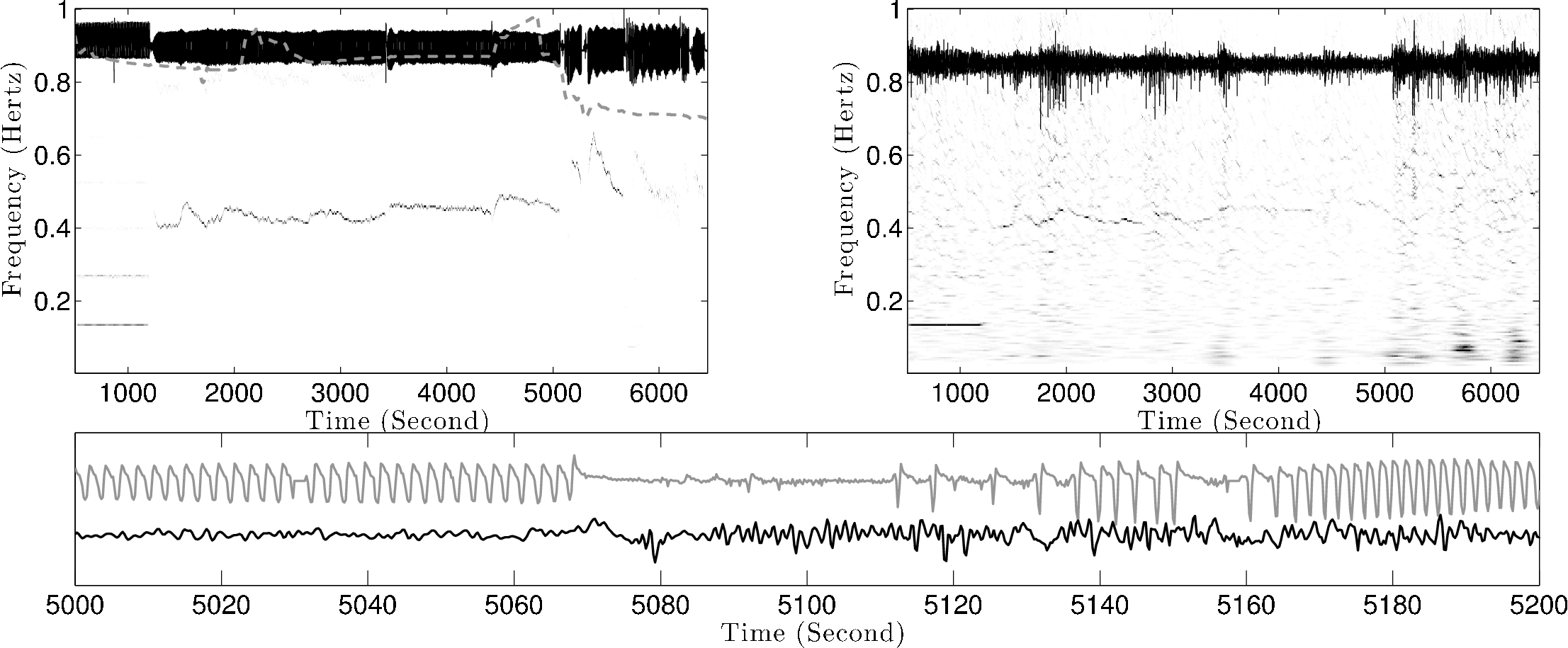}} 
\end{center}
\caption{Upper Left: the SST of the recorded respiratory signal, $Y_{\textup{Cflow}}$, from a subject with AF with $Y_{\textup{Cflow}}$ superimposed as the black curve and the end-tidal concentrations of anaesthetic gas superimposed as the gray dashed curve.  
The mechanical ventilator was removed at the 1205-th second and the subject started to breath spontaneously. The skin suture started right after removing the mechanical ventilator. The anesthetic gas was increased at the 2070-th second due to the hypertension induced by the skin suturing. The whole surgical intervention was finished at the 4870-th second and the anesthesiologist started to reduce the anesthetic medication concentration.
At the 5158-th second, the first reaction (movement of his head and neck) occurred. About 1 minutes before the first reaction, the respiratory pattern became irregular. 
Notice the dominant curve inside the band around $0.4$ Hertz before the $5158$-th second is corresponding to the IRR. The signal after the first reaction, however, does not always satisfy the model (\ref{eq:resp_model}), so we can see an interlacing appearance of a dominant curve, and the TF representation is blurred in between. Upper Right: the SST of the EDR signal, $Y_{\textup{EDR}}$, superimposed with $Y_{\textup{EDR}}$ as the black curve. We can see a dominant curve around $0.4Hz$ from $1200$-th second to $5100$-th second, which disappear afterward. This finding coincides with that in the SST of $Y_{\textup{Cflow}}$. 
Lower: the zoomed in signals around the first reaction (at the 5158-th second). The $Y_{\textup{Cflow}}$ is shown in the gray curve and $Y_{\textup{EDR}}$ is shown in the black curve. There is a re-calibration in $Y_{\textup{Cflow}}$ at the 5035-th second. Note that we can easily see visually the transition from the ``regular oscillation'' to the ''irregular oscillation'' in $Y_{\textup{Cflow}}$, while it is not easy to infer too much from $Y_{\textup{EDR}}$.}
\label{fig:12}
\end{figure}

Then, the analysis result of an anesthetized patient with AF under anesthesia is illustrated in Figure \ref{fig:12}. The SST of the recorded respiratory flow $Y_{\textup{Cflow}}$ is shown in the upper left sub-figure of Figure \ref{fig:12} with $Y_{\textup{Cflow}}$ superimposed as the black curve and the end-tidal concentrations of anaesthetic gas superimposed as the red curve. The SST of the $Y_{\textup{EDR}}$ is shown in the upper right subfigure of Figure \ref{fig:12} with $Y_{\textup{EDR}}$ superimposed as the black curve. 
Note that the transition from the regular pattern to non-regular pattern can be seen clearly both in$Y_{\textup{Cflow}}$ and its SST. On the other hand, this pattern transition can not be visually seen easily in $Y_{\textup{EDR}}$l shown in the lower subfigure of Figure \ref{fig:12}, while we can still see the transition in its SST figure. Indeed, before the first reaction at the 5158-th second, a visually viewable dominant curve exists in the TF representation, and the whole TF representation becomes blurred afterward. However, it is not an easy task to infer too much from $Y_{\textup{EDR}}$ as shown in the lower subfigure of Figure \ref{fig:12}.

In summary, the breathing dynamics of the respiration revealed from SSTEDR provides clinical information regarding the anesthetic depth, even for patients with AF. To further study this positive result and make it clinically helpful by quantifying the findings, a large scale study is ongoing and will be reported in our future work.

\section{Acknowledgements}
Hau-Tieng Wu acknowledges support by AFOSR grant FA9550-09-1-0643 and the valuable discussions with Professor Charles K. Chui and Professor Ingrid Daubechies.

\bibliographystyle{amsplain} 
\bibliography{EDR}

\clearpage

\end{document}